\documentstyle[twoside,11pt,amssymb]{article}
\title{\huge An Interpolation between Homology and\\
Stable Homotopy}

\author{Sadok Kallel\thanks{The author holds a Postdoctoral fellowship with 
PIMS and the University of British Columbia.}}
\date{}

\begin{document}
\maketitle


\def\la#1{\hbox to #1pc{\leftarrowfill}}
\def\ra#1{\hbox to #1pc{\rightarrowfill}}
\def\fract#1#2{\raise4pt\hbox{$ #1 \atop #2 $}}
\def\decdnar#1{\phantom{\hbox{$\scriptstyle{#1}$}}
\left\downarrow\vbox{\vskip15pt\hbox{$\scriptstyle{#1}$}}\right.}
\def\decupar#1{\phantom{\hbox{$\scriptstyle{#1}$}}
\left\uparrow\vbox{\vskip15pt\hbox{$\scriptstyle{#1}$}}\right.}
\def\ldnar#1{\phantom{\hbox{$\scriptstyle{#1}$}}
\left\downarrow\vbox{\vskip15pt\hbox{$\scriptstyle{#1}$}}\left.}
\def\updown#1{\phantom{\hbox{$\scriptstyle{#1}$}}
\left\updownarrow\vbox{\vskip15pt\hbox{$\scriptstyle{#1}$}}\right.}

\parskip=1pc
\def\za{\vrule height6pt width4pt depth1pt}
\def\lrar{{\ra 2}}

\def\sp#1{SP^{#1}}
\def\spy{SP^{\infty}}
\def\map#1{\hbox{Map}_{#1}}
\def\ext{\Lambda}

\def\bbc{{\Bbb C}}
\def\bbp{{\Bbb P}}
\def\bba{{\Bbb A}}
\def\bbz{{\Bbb Z}} 
\def\bbr{{\Bbb R}} 


\begin{abstract}

By considering labeled configurations of ``bounded multiplicity'', one
can construct a functor that fits between homology and stable
homotopy. Based on previous work, we are able to give an equivalent
description of this labeled construction in terms of loop space
functors and symmetric products. This yields a direct generalization
of the May-Milgram model for iterated loop spaces, and answers
questions of Carlsson and Milgram posed in the handbook.  We give a
classifying space formulation of our results hence extending an older
result of Segal.  We finally relate our labeled construction to a
theory of Lesh and give a generalization of a well-known theorem of
Quillen, Barratt and Priddy.

\end{abstract}


\noindent{\bf\Large \S1 Introduction}

In their paper in the {\sl Handbook of algebraic topology},
G. Carlsson and R.J. Milgram ([CM], \S 7) consider the following
space: Fix $k\geq 1$, $d\geq 1$ and let $F^d(\bbr^k,n)$ be the space
of ordered $n$-tuples of vectors in $\bbr^k$ so that no vector occurs
more than $d$ times in the $n$-tuple.  When $d=1$, $F^1(\bbr^k,n)$ is
the usual ordered configuration space $F(\bbr^k,n)$ consisting of
vectors of $\bbr^k$ with disjoint coordinates (cf. \S2).

Observe that the space $F^d(\bbr^k,n)\subset (\bbr^k)^n$, for a fixed $n\geq
1$, admits an action of the $n$-th symmetric group $\Sigma_n$ given by
permuting coordinates. Let $X$ be a topological space with a chosen
basepoint $*$.  We can then associate to the collection ${\cal F} =
\{F^d(\bbr^k,n)\}_{n\geq 0}$ a ``labeled'' construction
$$C^d(\bbr^k, X) = \coprod_{n\geq
0}F^d(\bbr^k,n)\times_{\Sigma_n}X^n/\sim \leqno{1.1}$$ 
where $\sim$ is a (standard) basepoint identification described in
\S2, and where $\times_{\Sigma_n}$ denotes the orbit space under the
permutation actions of $\Sigma_n$. (We make the convention here that
the term corresponding to $n=0$ is basepoint.)
The space $C^d(\bbr^k, X)$ is
(path) connected whenever $X$ is.

The authors in [CM] now raise the question of determining the homotopy
type of $C^d(\bbr^k, X)$ for connected $X$ and $d>1$. The case $d=1$
having been long known (see 1.3 below), an earlier attempt to answer
their question for the case $d=2$ was carried out by Karagueuzian
[Kr].

Let $\sp{d}(-)$ be the $d$-th symmetric product functor. This we
recall is defined as the quotient $\sp{d}(X) = X^d/\Sigma_d$ where
$\Sigma_d$ acts by permuting coordinates.  Our first main result
takes the form

\noindent{\bf Theorem} 1.2:~{\sl Let $X$ be (path) connected and CW.
Then there are homotopy equivalences for all $d\geq 1$;
$${C}^d(\bbr^k; X)\fract{\simeq}{\ra 3}\Omega^k(\sp{d}(\Sigma^k X)).$$}

\noindent{\bf Corollary} 1.3 (May-Milgram): {\sl 
Suppose $X$ path connected and CW. Then 
$C^1(\bbr^k, X)\simeq\Omega^k\Sigma^kX$, 
where $\Omega^k\Sigma^k(X)$ is the iterated $k$-fold loop space on
the $k$-fold suspension of the space $X$.}

\noindent{\bf Corollary} 1.4: In the case $d=k=1$, the space $C(\bbr;
X)$ ($X$ connected) is homotopy equivalent to the James construction
$J(X)$ described as the free monoid on points of $X$ with $*$ as the zero
element (see \S6). The equivalence 1.2 in this case yields James'
theorem
$$J(X)\simeq \Omega\Sigma X.$$

\noindent{\bf Remark} 1.5: It must be pointed out that an analogue of
theorem 1.2 based on little cubes of Boardman-Vogt has been obtained
independently by Fumiko Kato in her Master's thesis with D. Tamaki.
Our result also appears to have been known to F. Cohen and
C.F. Bodigheimer (although no proof of it has been published.)

\noindent{\bf Remark} 1.6:
Of course the space $\bbr^k$ in the construction 1.1. can be replaced
by any ``ground'' manifold $M$. When $M=*$ is the one point space,
$\Sigma_n$ acting trivially, the labeled space
$$\coprod_n *\times_{\Sigma_n}X^n/_{\sim}$$
is identified with the infinite symmetric product $\spy (X)$ (see
\S2).  It is well-known ([DT]) that $\pi_*(\spy (X))\cong {\tilde
H}_*(X;\bbz )$.  On the other hand and when $M=\bbr^{\infty}$, theorem
1.3 shows that
$C^1(\bbr^{\infty},X)\simeq\Omega^{\infty}\Sigma^{\infty}(X)=Q(X)$ and
hence that $\pi_*(C^1(\bbr^{\infty},X))\cong\pi_*^s(X)$. From this
perspective, the labeled constructs $C^d(M,-)$ for $d>1$ turn out to
provide intermediate functors between stable homotopy and integral
homology. This fact, which also explains our choice of title, provides
an extra incentive for the study of these functors.  

In general, replace $\bbr^k$ in 1.1 by an open manifold $M$ or by a
compact $M$ with non-empty boundary.  The following theorem extends a
result of Bodigheimer ([Bo])

\noindent{\bf Theorem} 1.7:~{\sl Let $M$ be a $k$-dimensional
manifold, smooth and stably parallelizable, with boundary $\partial
M\neq\emptyset$.  Assume $X$ (path) connected. Then there is a
homotopy equivalence
$${C}^d(M; X)\fract{\simeq}{\ra 3}
\map{}(M,\partial M; \sp{d}(\Sigma^k X))$$
where the mapping space on the right corresponds to all maps of $M$
into $\sp{d}(\Sigma^k X)$ sending $\partial M$ to basepoint.}

We next extend theorem 1.2 to the case when $X$ is disconnected. 
We observe (after Segal) that ${C}^d(\bbr^k; X)$ is homotopy
equivalent to an associative topological monoid ${\bar C}^d(\bbr^k;
X)$.  This monoid admits a classifying space $B{\bar C}$ and we
show (compare [S2] for the case $d=1$)

\noindent{\bf Theorem} 1.8:~{\sl Let $X$ be a topological space.
Then there is a weak homotopy equivalence
$$B{\bar C}^d(\bbr^k; X)\fract{\simeq_w}{\ra 3} 
\Omega^{k-1}(\sp{d}(\Sigma^k X)).$$}

\noindent{\bf Remark} 1.9: The case $X=S^0$ has been studied earlier
in [K] and in [GKY]. One writes $C^d(\bbr^k,S^0)=\coprod
C^d(\bbr^k,i)$ and this is a disconnected topological monoid (up to
homotopy) and hence we group complete it as follows. For $i\in\bbz^+$,
we define the space $\bbr^k_i$ to be the subset of $\bbr^k$ given as
follows
$$\bbr^k_i = \{(x_1,\ldots, x_k )\in\bbr^k,~|~0< x_k< i\}.$$
Clearly $F^d(\bbr^k,n)\cong F^d(\bbr^k_i,n)$ and for each $i$ we have 
an inclusion
$$\iota_i: F^d(\bbr^k_i,n)\lrar F^d(\bbr^k_{i+1},n+1), \zeta\mapsto \zeta+z_i$$
for some chosen $z_i\in (i,i+1)$ (here $+$ refers to juxtaposition).
The direct limit of the $\iota_i$ is denoted by
$F^d_{\infty}(\bbr^k)$. Of course this space admits an action of
$\Sigma_{\infty}$ and the quotient we denote by
$C^d_{\infty}(\bbr^{k})$. 
 
\noindent{\bf Theorem} 1.10: [K], [GKY]~{\sl Assume $d> 1$. Then there
is a homotopy equivalence
$${C}^d_{\infty}(\bbr^k)\simeq \Omega^k_0\sp{d}(S^k)$$
where $\Omega^k_0$ refers to the component of degree $0$ maps.
When $d=1$, one has instead a homology equivalence (a classical
result of Segal)
$$H_*({C}^1_{\infty}(\bbr^k);\bbz )\cong 
H_*(\Omega^k_0S^k;\bbz )$$}

To see how one can recover both 1.2 and 1.10 from theorem 1.8,
assume first that $X$ is connected. In this case, ${\bar C}= {\bar
C}^d(\bbr^k; X)$ is also connected and hence
$\Omega B{\bar C}\simeq {\bar C}$. This then yields theorem 1.2.
Theorem 1.10 on the other hand is deduced from 1.8 by letting $X=S^0$
and using the group completion theorem (cf. [KP], [MS]). 

In proving the theorems above we avoid entirely the theory of operads
(cf. [Ma]), iterated loop spaces (cf. [CM]) or classifying spaces as
in [S2]. Instead we use an interesting shortcut construction given in
the form of the ``scanning'' map. Such a construction has had much
interesting use recently (cf. [K], [GKY], [S1]) and one aim of
this paper is to advertize some of the ideas there.

The second part of this paper ties the above results
to recent work of K. Lesh [L1-2]. First notice that when $d=1$, then
${C}_{\infty}(\bbr^{\infty})$ is $B\Sigma_{\infty}$ up to homotopy.
The space $\coprod_{k\geq 0} B\Sigma_k$ has the structure of a
disconnected monoid with composition being induced from the pairings
$\Sigma_n\times\Sigma_m\rightarrow\Sigma_{n+m}$.  The
Baratt-Priddy-Quillen (BPQ) theorem states that
$$B\left(\coprod_{k\geq 1}B\Sigma_k\right)\simeq\Omega^{\infty-1}S^{\infty}
\leqno{1.11}$$
Alternatively, one recalls that $\Omega B\left(\coprod
B\Sigma_k\right)$ is the ``group completion'' $Gp$ of $\coprod_{k\geq
0} B\Sigma_k$. The B-P-Q theorem now states that $Gp\left(\coprod
B\Sigma_k\right)\simeq QS^0$ and if $Gp_0$ is anyone of the components
of $Gp$ (indexed by $\bbz$), then there are isomorphisms
$$H_*(Gp_0)\cong H_*(QS^0_0)\cong 
H_*(B\Sigma_{\infty})= H_*(C_{\infty}(\bbr^{\infty}))$$
where $QS^0_0$ is any component of $\Omega^{\infty}S^{\infty}$. 
From this point of view, the space $C_{\infty}(\bbr^{\infty})$ is closely
associated with the group completion of the ``family'' of groups
$\Sigma_k, k\geq 1$. Is there an analog of 1.11 for 
$C_{\infty}^d(\bbr^{\infty})$, $d>1$? Equivalently, which group
completion does the space $C_{\infty}^d(\bbr^{\infty}), d>1$ describe?

In \S5, we introduce the notion of a ``family'' of subgroups and of
their classifying spaces (a construction due to T. tom Dieck). It turns
out that given ``compatible'' families ${\cal F}_n$ (each consisting of
a collection of subgroups of $\Sigma_n$), one can associate
to them a topological monoid $\coprod B{\cal F}_n$ of which group
completion is an infinite loop space ([L1]). As a special case, we
consider the collection of subgroups
$$H_{i_1,\ldots, i_k
}=\Sigma_{i_1}\times\Sigma_{i_2}\times\cdots\times \Sigma_{i_k},~~
i_1+i_2+\cdots + i_k=n$$ 
of $\Sigma_n$.  We define a ``family'' ${\cal F}_n^d$ to consist of
all $H_{i_1,\ldots, i_k }$ with $i_j\leq d$ and $i_1+i_2+\cdots +
i_k=n$, together with their subgroups.  For every $n\geq 1$, the
family ${\cal F}^d_n$ affords a classifying space construction $E{\cal
F}_n^d$ (that is there is a $\Sigma_n$ space $E{\cal F}_n^d$ such that
the fixed point set under the action of $H_{i_1,\ldots, i_k }\in{\cal
F}_n^d$ is contractible, and otherwise it is empty).  The quotient
spaces $B{\cal F}_n^d = E{\cal F}_n^d/\Sigma_n$, $n\geq 1$, are
compatible in the sense that the disjoint union $\coprod B{\cal
F}_n^d$ has the structure of a topological monoid. We now prove
(cf. \S5)

\noindent{\bf Proposition} 1.12:~{\sl For all $d\geq 1$, we have the
following equivalence
$$B\left(\coprod_{n\geq 1}B{\cal F}_n^d\right)\simeq
\Omega^{\infty-1}\sp{d}(S^{\infty}).$$}

Again one may reinterpret this result by writing that $H_*(Gp_0)\cong
H_*(C^d_{\infty}(\bbr^{\infty}))$, where $Gp_0$ is any component of
$\Omega B (\coprod_{n\geq 1}B{\cal F}_n^d)$.  When $d=1$, ${\cal F}_n$
is the trivial family (consisting of the trivial subgroup in
$\Sigma_n$), $B{\cal F}_n=B\Sigma_n$ and one recovers 1.11 this way.
We note that theorem 1.12 is quite closely related to proposition 7.4 of [L2].

\noindent{\sc Acknowledgments:} The family of subgroups ${\cal F}_n^d$
mentioned above and defined in \S5 was suggested to the author by
Kathryn Lesh and the results of that section follow essentially from a
discussion with her. The author is grateful to her for discussing her
work with him and for commenting on an early version of this paper.
We also thank Dai Tamaki and in particular Denis Sjerve for his support
and for many stimulating discussions.


\vskip 10pt
\noindent{\bf\Large \S2. Labeled Constructions in Topology}

Constructions of the form 
$$\coprod_{n\geq 1} F(n)\times_{\Sigma_n}X^n/\sim\leqno{2.1}$$
where $F(n)$ is a space admitting a $\Sigma_n$ action and $\sim$ some
suitable identification, are called {\sl labeled constructions}. The
space $X$ is the {\sl label} space.  When $F(n)$ is some subspace of
$M^n$ where $M$ is some ``ground" space, 2.1
is referred to as a labeled {\sl configuration} space construction.
The first constructions of the sort trace their origin to the work of
Milgram on iterated loop spaces ([Mi]).

\noindent{\bf Notation and Terminology}: Traditionally, a
configuration of points of $M$ is any finite set of points of $M$
(that is {\sl unordered} and {\sl disjoint} points). However it is
better to enlarge this terminology to include symmetric products as
well; that is unordered points $\langle x_1,\ldots, x_n\rangle$,
$x_i\in M$, with the $x_i$ not necessarily distinct. We write $\langle
x_1,\ldots, x_n\rangle$ or $\sum n_ix_i$ to describe a configuration.

\noindent{\bf Example} 2.2: The ``Fadell'' space:~
$F(M,n ) =\{(x_1,\ldots, x_n)~|~x_i\in M, x_i\neq x_j, i\neq j\}.$

The first basic labeled (configuration) space construction one
encounters in topology is the infinite symmetric product construction.
Let $X$ be any space with a basepoint 
$x_0\in X$.  We denote by $\Sigma_n$ the $n$-th symmetric
group. Finally let $*$ be the one point space. Then we can consider
the quotient space
$$\spy (X) = \coprod_{k\geq 1}*\times_{\Sigma_k}X^k/\sim\leqno{}$$
where $\Sigma_k$ acts on $*$ trivially and on $X^k$ by permuting
coordinates, and where $\sim$ is the (basepoint) identification
$$*\times_{\Sigma_k}(x_1,\ldots, x_k)\sim *\times_{\Sigma_{k-1}}
(x_1,\ldots, {\hat x}_i,\ldots, x_k),~~\hbox{whenever}~~x_i=x_0.$$
Here $\hat x_i$ means deleting the $i$-th entry.

Similarly, it turns out that stable homotopy also admits a labeled
configuration space construction. This is given by the model
$$Q(X) = \coprod_{n\geq 1}F(\bbr^{\infty},n)\times_{\Sigma_n}X^n/\sim$$
where $F(\bbr^{\infty},n)\subset (\bbr^{\infty})^n$ is as defined in 2.2
and where $\sim$ is the basepoint relation
$$\small\displaystyle
(m_1,\ldots, m_n)\times_{\Sigma_n}(x_1,\ldots, x_n)\sim
(m_1,\ldots,{\hat m}_i,\ldots, m_n)\times_{\Sigma_{n-1}}
(x_1,\ldots,{\hat x}_i,\ldots, x_n)~~\hbox{if}~~
x_i=x_0.\leqno{2.4}$$
A stable version of the theorem of May and Milgram described in 1.3.
gives that $Q(X)\simeq \Omega^{\infty}\Sigma^{\infty}X$ for connected $X$
and this of course implies that $\pi_i(Q(X))\cong \pi_i^s (X).$

\noindent{\bf Remark} 2.5: The space $F(\bbr^{\infty},n)$ above is
characterized by the fact that it is contractible and admitting a {\sl
free} $\Sigma_n$ action. One can choose $F(\bbr^{\infty},n)$ to be a
model for $E\Sigma_n$ the acyclic bar construction on $\Sigma_n$ and
hence the quotient space
$$C(\bbr^{\infty},n)=F(\bbr^{\infty},n)/\Sigma_n$$
is identified with the classifying space $B\Sigma_n$.
Naturally the constant map $F(\bbr^{\infty},n)\rightarrow *$ induces an
equivariant map $F(\bbr^{\infty},n)\times X^n\rightarrow *\times X^n$
and hence a map $QX\lrar\spy (X)$
which turn out to be ``a space level representation''
of the Hurewitz map
$$\pi_i^s(X)\lrar H_i(X).$$

\noindent{\bf Remark} 2.6: One notices that singular homology and
stable homotopy sit at opposite ends of this labeled configuration
space construction. Indeed, the basepoint $*$ which one associates to
integral homology, is a space which admits the ``most unfree" action
of $\Sigma_n$ possible. On the other hand, $F(\bbr^{\infty},n)$ which
is associated to $\pi_*^s$, is a space on which $\Sigma_n$ acts
freely.  Naturally one can consider ``in between'' spaces $F(n)$;
$$*\subset F(n)\subset F(\bbr^{\infty},n)$$ 
which admit an action of $\Sigma_n$ and obtain in this way some
interpolations between both theories. This paper deals with one
particular such interpolation.


\vskip 10pt
\noindent{\bf\Large\S3 The functor $C^d(-;X)$ and Its Completion}

We assume $M$ to be a connected manifold (either open or
compact with non-empty boundary.)  Let $F^{d}(M,n)\subset M^n$ be the
set of $n$-tuples of points of $M$ such that no point occurs more than
$d$ times in the tuple;
\begin{eqnarray*}
F^d(M,n)=\bigcup~\{(\underbrace{x_1,\ldots, x_1}_{i_1}, 
\underbrace{x_2,\ldots, x_2}_{i_2},\cdots, 
\underbrace{x_k,\ldots, x_k}_{i_k})&|&x_i\in M, x_i\neq x_j, i\neq j,\\
&&~\hbox{and}~i_j\leq d,~~i_1+\cdots +i_k = n \}.
\end{eqnarray*}
When $d=1$, $F^1(M,n)$ is just the space $F(M,n)$
described in 2.2 and one has the filtration
$$F(M,n)=F^1(M,n)\subset F^2(M,n)\subset\cdots\subset F^n(M,n)=M^n.$$

\noindent{\sc Example}: Let $M=\bbr^2$. Then $C^d(\bbr^2,n)$
consists of $n$ points in the plane each having multiplicity less than
$d$.  Each such configuration is uniquely identified with the roots of
a monic complex polynomial and so $C^d(\bbr^2,n)$ is identified with
the set of degree $n$ complex (monic) polynomials having roots of
multiplicity less than $d$; a space originally studied by Vassiliev.

As before, the family $\{F^d(M,n)\}_{n\geq 0}$ admits a labeled space
construction
$$C^d (M;X)=\coprod_{n} F^d(M,n)\times_{\Sigma_n}X^n/_{\sim}\leqno{3.1}$$
where $\sim$ is the basepoint relation described in 2.4. 
Notice that this labeled construction is a subset of $\spy (M\ltimes
X)$ where $M\ltimes X$ is the half-smash obtained from $M\times X$ by
collapsing $M\times *$. To see this inclusion, simply rewrite an
element $(m_1,\ldots, m_n)\times_{\Sigma_n}(x_1,\ldots, x_n)$, $m_i\in
M$, $x_i\in X$, in the form $\sum^n (m_i, x_i)$.

The labeled construction 3.1 defines a bifunctor $C^d(-;-)$ which is
a homotopy functor in $X$ and an isotopy functor in $M$ (that is $M$
must be in the category of spaces and injective maps).  Of course when
$X=S^0$, we write
$$C^d(M,S^0):= C^d(M) = \coprod_{k\geq 0}F^d(M,k).$$

We observe that we can put a monoid structure on $C^d(\bbr^k, X)$
up to homotopy. This is done as follows:
Let $\bbr^k_t$ be given as in \S1;
$$\bbr^k_t = \{(x_1,\ldots, x_n )\in\bbr^k,~|~0< x_n< t\}\leqno{3.2}$$
and define the space
$${\bar C}^d(\bbr^k;X) = \{(\zeta, t)\in C^d(\bbr^k_t;X)\times \bbr^+\}.
\leqno{3.3}$$
Again we have that ${\bar C}^d(\bbr^k,X)\simeq C^d(\bbr^k, X)$ and
this new modified space has now the structure of an associative
topological monoid with a composition law given by juxtaposition
$$ C(\bbr^k_t;X)\times C(\bbr^k_{t'};X)\lrar C(\bbr^k_{t+t'});X),~~
(\zeta,\zeta')\mapsto (\zeta + T_t\zeta')$$
where $T_t$ is translation $(0,t'))\lrar (t, t+t')$.
In the case $d=k=1$ this is the same as the well-known James
construction.

\noindent{\bf Lemma} 3.4 ([S2]):~{\sl Let $J(X)$ be the free monoid
generated by points of $X$ (with $*$ as the zero element). Then there
is a homomorphism ${\bar C}(\bbr,X)\lrar J(X)$ which is a homotopy
equivalence.}

\noindent{\sc Proof:} Here we have
$${\bar C}(\bbr,X) = \{(\zeta, t)\in C((0,t);X)\times \bbr^+\}$$
Let $(\zeta, t)\in {\bar C}(\bbr,X)$, where $\zeta$ is of the form
$$\zeta = (t_1,\ldots, t_n)\times_{\Sigma_n}(x_1,\ldots, x_n),~
\hbox{for some}~n, 0<t_i<t.\leqno{3.5}$$ 
We can then construct the map
$$\tau: {\bar C}(\bbr,X)\lrar J(X),~~(\zeta, t)\mapsto x_1+x_2+\cdots
+x_n.$$
This is well defined for one checks that the identifications are
preserved when $x_i=*$. One also checks that this is a
homomorphism. Indeed, the pairing $\mu$ on ${\bar C}(\bbr,X)$ is given
by
\begin{eqnarray*}
  \mu \left((\zeta, t), (\zeta', t')\right)&\mapsto& (\zeta+T_t\zeta',
  t+t')\\
  &&\hbox{where}~\zeta+T_t\zeta' = (t_1,\ldots, t_n,
  t+t'_1,\ldots, t+t'_m)\times_{\Sigma_{n+m}} (x_1,\ldots, x_n,
  x'_1,\ldots, x'_m).
\end{eqnarray*} 
It follows that $\tau\mu \left((\zeta, t), (\zeta', t')\right) =
x_1+\cdots + x_n + x'_1+\cdots + x'_m = \tau\mu (\zeta, t) + \tau\mu
(\zeta',t')$ as desired.

We finally need verify that $\tau$ is a homotopy equivalence.  Notice
that the point $\zeta$ in 3.5 above can be written as an ordered tuple
$$\zeta = ((t_1,x_1),\ldots, (t_n,x_n)),~\hbox{where}~
t_1<t_2<\cdots <t_n<t.$$
Consider the map 
$$\alpha: {\bar C}(\bbr,X)\lrar J(X)\times\bbr^+ ,~
(\zeta, t)\mapsto (\tau (\zeta, t), t).$$
Then the preimage of a point $(x_1+\ldots +x_n, t)$
under $\alpha$ is an open simplex
$$\Delta_n=\{(t_1,\ldots, t_n), 0<t_1<t_2<\cdots <t_n<t\}$$ 
and this is contractible. The signifies that $\alpha$ is acyclic and
hence a homotopy equivalence. It follows that $\tau$ is an equivalence
and the proof is complete.  \hfill\za

\noindent{\bf Group Completion}: Since ${\bar C}^d = {\bar
C}^d(\bbr^k;X)$ is a monoid (possibly disconnected), it admits a group
completion $\Omega B{\bar C}^d$. A very handy way of describing this
group completion in terms of configuration spaces is given as follows.

We shall suppose that $X$ has finitely many components. This implies
that ${\bar C}^d(\bbr^k;X)$ has $\bbz^m$ components, with
$m=|\pi_0(X)|-1$.  For each $i\leq |\pi_0(X)|$ choose a point $p_i$ in
the $i$-th component of $X$ (you can let basepoint be the point in the
component of the identity). Let $\bbr^k_t$ be as in 3.2 and choose a
point $z_i\in\bbr^k_{i,i+1} =\bbr^k_{i+1}-\bbr^k_i$,
for all $i\geq 1$. We can then consider
the inclusion
$$\matrix{
{\bar C}^d(\bbr^k;X)&\fract{\tau_i}{\ra 2}&{\bar C}^d(\bbr^k;X)\cr
(\sum (m_r,x_r), i)&\mapsto& (\sum (m_r, x_r) + (z_i, p_i), i+1)
\cr}\leqno{3.6}$$
where $\sum (m_r, x_r)\in C^d(\bbr^k_i;X)$ is as described in 3.1.
The direct limit over these maps is denoted by 
${\hat C}^d(\bbr^k;X)$.

\noindent{\bf Lemma} 3.7:~{\sl Let $X$ be CW. Then $H_*({\hat
C}^d(\bbr^k;X))\cong H_*(\Omega B {\bar C}^d(\bbr^k;X)).$}

\noindent{\sc Proof}: If we let $\pi = H_0({\bar C}^d)$, then a
theorem of Kahn and Priddy states that
$$H_*({\bar C}^d)[\pi^{-1}]\cong H_*(\Omega B{\bar C}^d)\leqno{3.8}$$
where the left hand side means localization with respect to the
multiplicative set $\pi$.  The idea here of course is that by
inverting $\pi$, we are ``turning'' multiplication by elements of
$\pi$ into isomorphisms. That this is necessary is clear since $\Omega
B{\bar C}^d$ is a group and hence the image of $\pi$ under
$M\rightarrow \Omega B{\bar C}^d$ must consist of units.

Now notice that the point $(z_i, p_i)\in {\bar C}^d(\bbr^k,X)$ constructed
above represents a point $e_i\in\pi_0({\bar C}^d)$. The stabilization
maps in 3.6 correspond therefore to maps
$${\hat C}^d(\bbr^k;X)\simeq\lim_{\lrar\atop e_i\in\pi}\left(
{\bar C}^d(\bbr^k;X)\fract{\cdot e_i}{\ra 4}{\bar C}^d(\bbr^k;X)\right).$$
and this direct limit (by construction) 
must satisfy
$H_*({\hat C}^d(\bbr^k;X))\cong H_*({\bar C}^d(\bbr^k;X))[\pi^{-1}].$
The claim follows from 3.8.
\hfill\za


\vskip 10pt
\noindent{\bf\Large\S4 The Correspondence and Proofs of Main Theorems}

\noindent{\bf Scanning Labeled Configurations}: Details of this
construction can be found in [K]. Let $M$ be smooth of
dimension $k$ and suppose it is parallelizable. This means among other
things that every small (closed) neighborhood of $x\in M$ can be
canonically identified (via the exponential map for example) with a
(closed) disc $D^k$. Given a configuration $\zeta = \sum (m_i,x_i)\in
C^d(M;X)$, then its ``restriction'' to a neighborhood $D^k(x)\subset
M$ of $x$ gives rise to a new configuration 
$$\zeta\cap D^k(x)\in C^d(D^k(x); X).$$
The correspondence $\zeta\mapsto \zeta\cup D^k(x)$ is not continuous
as is, however by stipulating that whenever the points $m_i$ in 
$\zeta = \langle (m_1,x_1),\ldots, (m_n,x_n)\rangle$ get
close to $\partial D^k(x)$ they are ``dropped out''. We then
get a quotient map and a correspondence
$$\zeta\in C^d(M;X)\mapsto \zeta_x\in C^d(D^k(x),\partial D^k(x);X)$$
where ${C}^d(D^k(x),\partial D^k(x); X)$ is the identification space 
$$C^d(D^k,\partial D^k; X)=
\coprod_nF^d(D^k, n)\times_{\Sigma_n}X^n/_{\sim},~x\in X$$
and $\sim$ consist of the usual basepoint relation
together with the additional identification
$$(m_1,\ldots, m_n)\times_{\Sigma_n}(x_1,\ldots, x_n)
\sim (m_1,\ldots, {\hat m}_i,\ldots, m_n)
\times_{\Sigma_n}(x_1,\ldots, {\hat x}_i,\ldots, x_n)$$
whenever $m_i\in\partial D^k$ (i.e. when points of $D^k$ tend to the
boundary they are discarded together with their labels). 
Notice that $C^d(D^k,\partial D^k; X)\simeq C^d(S^k,*; X)$ and this
space has a canonical basepoint which we also write as $*$. 
Therefore and by scanning every neighborhood $D^k(X)$ of $x\in M$ and
identifying canonically the pair
$(D^k(x),\partial D^k(x))$ with $(D^k,\partial D^k) = (S^k,*)$ 
we construct a map 
$$S: C^d(M;X)\lrar\map{}(M, C^d(S^k,*;X)).$$
If in addition $M$ has a non-empty boundary $\partial M$, it is still
possible to scan labeled configurations
away from the boundary so that one obtains a map
$$S: C^d(M;X)\lrar\map{}(M/\partial M, C^d(S^k,*;X))$$
where the right hand side is the space of maps sending $\partial M$
to the basepoint in $C^d(S^k,*;X))$.

\noindent{\bf Lemma} 4.1:~{\sl There is a homotopy equivalence
$C^d(S^k,*; X)\simeq \sp{d}(\Sigma^k X).$}

\noindent{\sc Proof:} This is the analog of the unlabeled case
discussed in [K]. Choose $*$ in $S^k$ to be the south pole and let $x_0$
refer to the north pole. Let $U_{\epsilon}$ be an epsilon
neighborhood of $x_0\in S^k$. Notice that there is a radial homotopy,
injective on the interior of $U_{\epsilon}$ that expands the north
cap $U_{\epsilon}$ over the sphere and takes $\partial U_{\epsilon}$
to $*$.

Consider at this point the subspace $W_{\epsilon}$ consisting of
$(m_1, \ldots, m_n)\times_{\Sigma_n}(x_1\ldots , x_n)$ in
$C^d(S^k,*;X)$ such that at most $d$ points of $\langle m_1,\ldots,
m_n\rangle$ lie inside $U_{\epsilon}$.  By definition of
$C^d(S^k,*;X)$ each of its elements must fall into a $W_{\epsilon}$
for some $\epsilon$ and hence
$$C^d(S^k,*;X)= \bigcup_{\epsilon}W_{\epsilon}.$$
Now using the radial retraction described above, each labeled
configuration $(m_1, \ldots, m_n)\times_{\Sigma_n}(x_1\ldots , x_n)$
is sent to (at most) a $d$ point configuration which takes the form
$\sum_1^d (m_i,x_i)$. The identifications are such $(m,*)\sim *$ and
$(*,x)\sim *$ (here the $*$'s refer to the corresponding basepoints in
$X$, $S^k$ and $C^d$). The configuration $\sum_1^d (m_i,x_i)$ is
clearly an element of $\sp{d}(S^k\wedge X)=\sp{d}(\Sigma^k X)$ and the
lemma follows.  \hfill\za

Let $C^d_{\infty}(\bbr^k; X)$ denote the stabilized space described in
3.6. When $X$ is connected, one has the equivalence
$C^d_{\infty}(\bbr^k;X)\simeq C^d(\bbr^k;X)$.

\noindent{\bf Proposition} 4.3:~ {\sl For $X$ a topological space,
scanning induces an (integral) homology equivalence
$$S_*: H_*({\hat C}^d(\bbr^n; X))\fract{\cong}{\ra
3}H_*(\Omega^n\sp{d}(\Sigma^nX)).$$ 
This is a homotopy equivalence whenever $X$ is connected.}

\noindent{\sc Proof}: The proof uses properties of the bifunctor
$C^d(-,-)$ (cf. appendix) and is based on an induction on a suitable
handle decomposition of $D^n$ (the decomposition as described below is
due to Bodigheimer [Bo].)

Let $H^k=D^n=D^k\times D^{n-k}$ and write $A^k = S^{k-1}\times
D^{n-k}$ (note that $A^0=\emptyset$).  ($H^k$ is called a
handle. Handles make up manifolds by sequences of attachments. the
handle $H^k$ has ``index'' $k$ if $A^k$ is the part of its boundary
along which it is to be attached.)

Let $I_k\subset D^n=[0,1]^n$ denote the subset of $(y^1,\ldots, y^n)$
such that $y^i=0$ or $y^i=1$ for some $i=1,\ldots, k-1$, or $y^k=1$
(that is $I_k$ consist of all the boundary faces of $D^k \subset
D^n=D^k\times D^{n-k}$ safe the face $y^k=0$). Now let
$H_k=[0,1]^{k-1}\times [0,{1\over 2}]\times [0,1]^{n-k}$.  Then there
is a cofibration sequence
$$
(H_k, H_k\cap I_k)\lrar (D^n, I_k)\lrar (D^n, H_k\cup I_k).
$$
The pair $(H_k, H_k\cap I_k)$ can be identified with $(D^n,
S^{k-2}\times D^{n-k+1})$ hence representing a $k-1$-handle
$(H^{k-1},A^{k-1})$, while the pair $(D^n, H_k\cup I_k)=(D^n,
S^{k-1}\times D^{n-k})$ represents a handle $(H^k,A^k)$ of index $k$.
Applying the functor $C^d(-;X)$ and then scanning yields the
homotopy commutative diagram for all $k\geq 0$
$$\matrix{
   C^d(H^{k},A^{k};X)&\lrar&\Omega^{n-k}\sp{d}(\Sigma^nX)\cr
   \downarrow&&\downarrow\cr
   C^d (D^n,I_k;X)&\lrar& PS \cr
   \decdnar{\pi}&&\decdnar{}\cr
   C^d(H^{k+1},A^{k+1};X)&\lrar&\Omega^{n-k-1}\sp{d}(\Sigma^nX)\cr}$$
When $k>0$, $A^k\neq\emptyset$ and 
proposition A.3 asserts that the left hand side is a quasifibration.
In this case one shows (inductively) that
$C^d(H^k,A^k;X)\simeq \Omega^{n-k}\sp{d}(\Sigma^nX)$.
When $k=0$, $A^k=\emptyset$ and we need to pass to group completed spaces
${\hat C}^d$. For that we assume $X$
has finitely many components first. We then have the diagram
$$\matrix{
   {\hat C}^d(D^n;X)&\lrar&\Omega^{n}\sp{d}(\Sigma^nX)\cr
   \downarrow&&\downarrow\cr
   {\hat C}^d (D^n,*;X)&\fract{\simeq}{\lrar}& PS \cr
   \decdnar{\pi}&&\decdnar{}\cr
   C^d(H^1,A^{1};X)&\fract{\simeq}{\lrar}&\Omega^{n-1}
  \sp{d}(\Sigma^nX)\cr}$$
where the left hand side is now a homology fibration (by proposition A.5).
Since $\Omega^{n}\sp{d}(\Sigma^nX)$ must be the homotopy fiber of 
${\hat C}^d (D^n,*;X)\lrar C^d(H^1,A^{1};X)$, it follows that
$H_*({\hat C}^d(D^n;X))\cong H_*(\Omega^{n}\sp{d}(\Sigma^nX))$ as
asserted.

The case of a general disconnected $X$ is obtained once we observe
that $X$ can be described as the direct limit of spaces $X_i$ with
$\pi_0(X_i)$ finite. The space $C^d_{\infty}(D^n;X)$ is then the
direct limit of the $C^d_{\infty}(D^n;X_i)$ and since direct limits
commute with homology, the claim follows.  
\hfill\za

\noindent{\bf Remark} 4.4: Proposition 4.3 can be stated with
$\bbr^k$ replaced with any manifold with non-empty boundary or with
any open manifold. The proof of the more general statement uses a
handle decomposition argument of the manifold $M$ and inductive use
of A.5. This yields theorem 1.7 of the introduction.

\noindent{\bf Corollary} 4.5: {\sl For $X$ a topological space,
there is a weak homotopy equivalence
$$B{\bar C}^d(\bbr^k;X)\simeq \Omega^{k-1}\sp{d}(\Sigma^kX).$$}

\noindent{\sc Proof:} When $X$ is connected, $\Omega B{\bar C}^d\simeq
{\bar C}^d$ and the result follows directly from 4.3. Suppose that $X$
disconnected, hence so is ${\bar C}^d$ as a topological monoid.  
Assume for now that $\pi_0(X)$ finite and let
${\hat C}^d(\bbr^k,X)$ be as defined in \S3.
The following sequence of isomorphisms
$$H_*(\Omega B{\bar C}^d)\fract{3.7}{\la 3}
H_*({\hat C}^d (\bbr^k;X))
\fract{4.3}{\ra 3} H_*(\Omega^k\sp{d}(\Sigma^kX))$$
shows that 
$H_*(\Omega B{\bar C}^d)\cong H_*(\Omega^k\sp{d}(\Sigma^kX))$.
This implies that
$$H_*(B{\bar C}^d)\cong H_*(\Omega^{k-1}\sp{d}(\Sigma^kX)).$$
Oberve that the fundamental groups map isomorphically
\begin{eqnarray*}
\pi_1(B{\bar C}^d)\cong \pi_0({\bar C}^d)&\cong&
\bbz^{|\pi_0(X)|-1}\cong 
{\tilde H}_0(X) = H_k(\Sigma^kX)\cong H_k(\sp{d}(\Sigma^kX))\\
&\cong&\pi_k(\sp{d}(\Sigma^kX)) = \pi_1\left(
\Omega^{k-1}\sp{d}(\Sigma^kX)\right)
\end{eqnarray*}
and the claim follows from Whitehead's criterion.

To treat the general case of $\pi_0(X)$ not necessarily finite, we
write again ${\bar C}^d_{\infty}(\bbr^k;X)$ as a direct limit
of spaces ${\bar C}^d_{\infty}(\bbr^k;X_i)$ for some $X_i$ with finite
$\pi_0$. Since homotopy groups commute with direct limits, the 
theorem follows.
\hfill\za


\vskip 10pt
\noindent{\bf\Large \S5 A Theory of Infinite Loop Spaces}

In this section we describe a construction of K. Lesh (cf. [L1-2])
which associates to a (compatible) family of groups an infinite loop
space. We then describe how our labeled construction fits in and
prove proposition 1.12 of the introduction.

Let $G$ be a group and let $\cal F$ be a collection of subgroups of
$G$ which is closed under conjugation; meaning that \hfill\break
$\bullet$ If $H\in {\cal F}$ and $g\in G$, then $g^{-1}Hg\in {\cal F}$
\hfill\break $\bullet$ If $H\in{\cal F}$ and $K$ a subgroup of $H$,
then $K\in {\cal F}$.\hfill\break Such a collection is called a {\sl
family}.

The prototypical example of a family would be to take all subgroups of
a group $G$ (a variant is to consider only the finite subgroups).  A
less trivial example would be to consider the family of elementary
abelian $p$-subgroups of $\Sigma_n$ which are generated by disjoint
$p$-cycles together with their subgroups. This family is studied
extensively in [L1].

It turns out that to a family $\cal F$ of subgroups of a group $G$
there is associated a {\sl classifying} space $B{\cal F}$ by work of
T. tom Dieck. More precisely, tom Dieck constructs a $G$-space $E{\cal F}$
with the property that the fixed point set $E{\cal F}^H$ of $H\subset
G$ is such that
$$
E{\cal F}^H\simeq *~~\hbox{for}~H\in{\cal F},~~~\hbox{and}~~~
E{\cal F}^H=\emptyset~~\hbox{for}~H\not\in{\cal F}
$$
Note that $E{\cal F}$ is always contractible since $*\in {\cal F}$ for
any family.  Naturally one then defines the classifying space $B{\cal
F}$ to be the orbit space of the $G$ action on $E{\cal F}$.

\noindent{\sc Example} 5.1: Let $\cal F$ consists only of the trivial
subgroup in $G$. Then $E{\cal F}=EG$.

We now specialize to the symmetric groups $\Sigma_n$ and we suppose
that for each $n$ we are given a family ${\cal F}_n$ of subgroups for
$G=\Sigma_n$. We recall that given two subgroups $H\in\Sigma_n$ and
$K\in\Sigma_m$, we have a group $H\times K\in\Sigma_{n+m}$ obtained as
the image of the composite
$$H\times K\hookrightarrow \Sigma_n\times\Sigma_m\hookrightarrow
\Sigma_{n+m}.$$

\noindent{\sc Definition} 5.2: The families $\{{\cal F}_n\}_{n\in{\bf
Z}^+}$ are {\sl compatible} if whenever $H\in {\cal F}_n$ and
$K\in{\cal F}_m$, then $H\times K\in{\cal F}_{n+m}$.

\noindent{\bf Theorem} 5.3 (Lesh):~{\sl Let $\{{\cal F}_n\}_{n\in{\bf
Z}^+}$ be a compatible choice of families, then $\coprod B{\cal F}_n$
has a monoid structure whose group completion is an infinite loop
space $L{\cal F}$. Such a space comes equipped with (natural) maps
$$QS^0\lrar L{\cal F}\lrar \bbz.$$}

\noindent{\sc Example} 5.4: Let ${\cal F}_n$ be the family consisting
of the trivial subgroup in $\Sigma_n$. Then $B{\cal F}_n=B\Sigma_n$
and so $L{\cal F}$ in this case is the group completion of $\coprod
B\Sigma_n$ which is known to correspond by a theorem of Barratt-Priddy
and Quillen to the infinite loop space $QS^0$.

We now relate the above constructions to the spaces
$C^d(\bbr^{\infty},n)$ and their stable version $C^d(\bbr^{\infty})$
constructed in \S3. Given $n\geq 1$ we consider the following
subgroups of $\Sigma_n$;
$$H_{i_1,\ldots, i_k}
=\Sigma_{i_1}\times\Sigma_{i_2}\times\cdots\times\Sigma_{i_k}\subset\Sigma_n,
~~i_j\leq n~\hbox{and}~i_1+i_2+\cdots + i_k=n.$$ Each such subgroup
$H_{i_1,\ldots, i_k}$ acts on $F^d(\bbr^{\infty},n)$ by permuting
points. Let
$${\cal F}_n^d = \{H_{i_1,\ldots, i_k}~~|~~i_j\leq d,~~
i_1+i_2+\cdots + i_k= n,~~\hbox{together with
their subgroups}\}.$$
It is not hard to see that ${\cal F}_n^d$ satisfies the conditions of
a family, and that the newly obtained 
families $\{{\cal F}_n^d\}_{n\in{\bf Z}^+}$ 
form a compatible collection.

\noindent{\bf Lemma} 5.5:~{\sl $E{\cal F}_n^d\simeq F^d(\bbr^{\infty},n)$.}

\noindent{\sc Proof:}
Pick $H= H_{i_1,\ldots, i_k}\in{\cal F}_n^d$ ($i_j\leq d$
and $\sum i_j=n$.) Then 
$$(E{\cal F}_n^d)^{H_{i_1,\ldots, i_k}} = 
\{(\underbrace{x_1,\ldots, x_1}_{i_1}, 
\underbrace{x_2,\ldots, x_2}_{i_2},\cdots, 
\underbrace{x_k,\ldots, x_k}_{i_k})~|~x_i\in \bbr^{\infty}\}$$
where the $x_i$ need not be distinct. If we let $X(\bbr^j)$ be the
subset of $(E{\cal F}_n^d)^H$ consisting of the
$x_i\in\bbr^j\subset\bbr^{\infty}$, then we see that $X(\bbr^j)$ is
open in $\bbr^{jn}$ and is the complement of hyperplanes of codimension 
at least $j$ (implying in particular that it is $j-2$
connected).  Since $(E{\cal F}_n^d)^H$ is the direct limit of
$X(\bbr^j)\hookrightarrow X(\bbr^{j+1})$ it must be contractible and
$(E{\cal F}_n^d)^H\simeq *$ as desired. What is left to show is that
the fixed point set of $H\not\in {\cal F}_n^d$ is empty. Observe that
any such $H$ must contain a cycle on (at least) $d+1$ letters. The
fixed points of such a cycle consists of configurations containing a
$d+1$ (or maybe more) repeated point.  Such a configuration cannot
exist in $F^d(\bbr^{\infty},n)$ (by definition) and $(E{\cal
F}^d)^H=\emptyset$.  \hfill\za

\noindent{\bf Proposition} 5.6:~{\sl For all $d\geq 1$, we have the
following equivalence
$$B\left(\coprod_{n\geq 1}B{\cal F}_n^d\right)\simeq
\Omega^{\infty-1}\sp{d}(S^{\infty}).$$}

\noindent{\sc Proof:} Theorem 1.8 still holds true for $k=\infty$
by a direct limit argument.  On the other hand, 5.5 implies that
$$\coprod_{n\geq 1}B{\cal F}_n^d 
=\coprod_{n\geq 1}C^d(\bbr^{\infty},n) = C^d(\bbr^{\infty};S^0)
\simeq {\bar C}^d(\bbr^{\infty};S^0)$$
and the claim follows.  \hfill\za


\vskip 10pt
\noindent{\bf\Large Appendix: Properties of the Bifunctor $C^d(-,-)$}

We here summarize the most useful properties of these functors.  The
details as well as a more extensive study of a more general class of
functors (so-called {\sl particle} functors) are given in [K].

\noindent{\bf Definition}: Let $f:X\lrar Y$ be a map of CW complexes.
For $x\in Y$, let $\iota_x$ denote the homotopy
inclusion of the preimage $f^{-1}(x)$ into the homotopy
fiber. Then\hfill\break $\bullet$ $f$ is a quasifibration if $\iota_x$
is a homotopy equivalence ($\forall x$)\hfill\break 
$\bullet$ $f$ is a homology equivalence if $\iota_x$ induces an
isomorphism in homology groups ($\forall x$).

Let now $M$ be an $n$-dimensional manifold, $M_0$ a closed submanifold
and $A\subset M$ a closed
subset. We assume that $N$ has a tubular neighborhood $U$ which
retracts to it (if it is a submanifold for example). We also assume
throughout that $M$ and $N$ are connected.
The following sequence of quotient maps 
$$(N, N\cap M_0)\lrar (M,M_0)\lrar (M,N\cup M_0)\leqno{A.1}$$
forms then a cofibration. Applying $C^d(-;X)$ to all terms of A.1
yields
$$C^d(N, N\cap M_0;X)\lrar C^d(M,M_0;X)\lrar C^d(M,N\cup M_0;X)
\leqno{A.2}$$
One of the main results about the functor $C^d(-)$ is that

\noindent{\bf Theorem} A.3: {\sl Let $N,M,M_0$ be as above.  If either
$X$ is connected or $N\cap M_0\neq\emptyset$, then the sequence A.2 is
a quasifibration.}

A short proof of this theorem can be found in [Bo] (for the case
$d=1$).  We are here concerned with the case when $N\cap
M_0=\emptyset$ or $X$ is disconnected. Let's assume $X=S^0$ for now
and let's consider the sequence
$$C^d(N)\lrar C^d(M,M_0)\lrar C^d(M,N\cup M_0).$$
This turns out not to have any particular structure. However by passing
to group completed space ${\hat C}^d$ as defined in \S3 we can show

\noindent{\bf Proposition} A.4: {\sl Let $M$ and $N$ be connected
manifolds, $N\subset M$, $M_0\subset M$ and $N\cap M_0=\emptyset$.
Assume $N$ has an end (or a boundary). Then
$${\hat C}^d(N)\lrar {\hat C}^d(M,M_0)\lrar C^d(M,N\cup M_0)$$
is a quasifibration if $d>1$ and a homology fibration if $d=1$}

\noindent{\sc Sketch of Proof:} 
The point here is that showing the sequence of spaces above
is a homology fibration (resp. a quasifibration) boils down to
showing that maps
$${\hat C}^d(N)\fract{+}{\ra 4} {\hat C}^d(N)$$
given by adjointing a given set of configurations is a homology
equivalence (resp. a weak homotopy equivalence).  Because of the very
construction of $\hat C$, adding configurations simply switches
components and since these components are the same, ``addition''
induces a homology equivalence. This is not (necessarily) a homotopy
equivalence since there is no obvious map backwards
(``subtraction'') which when composed with addition induces the
identity on components. When $d>1$ the fundamental group of
$C^d(N)$ abelianizes and this is enough to induce a homotopy
equivalence. For details the reader can consult [K].
\hfill\za

\noindent{\bf Proposition} A.5:~{\sl We let $M=\bbr^n$, $N, M_0\subset
M$ such that $M\cap N=\emptyset$. Let $X$ be CW with a finite number
of components.  Then the following is a homology fibration
$${\hat C}^d(M;X)\lrar {\hat C}^d(M, M_0;X) \lrar C^d(M,M_0\cup N;X).$$}


\addcontentsline{toc}{section}{Bibliography}
\bibliography{biblio}

\begin{thebibliography}{doC}

\bibitem{} [B] C. Berger, ``Op\'erades cellulaires et espaces de lacets
it\'er\'es'', Ann. Inst. Fourier, Grenoble {\bf 46} 4 (1996) 1125--1157.
\bibitem{}[Bo] C.F. Bodigheimer, ``Stable splittings of mapping spaces'', 
Algebraic topology, Proc. Seattle (1985), Springer lecture notes {\bf 1286}, 
174--187.
\bibitem{}[CM] G. Carlsson, R.J. Milgram, ``Stable homotopy and
iterated loop spaces'', Handbook of algebraic topology, Elsevier 1995.
\bibitem{} [D] A. Dold, ``Decomposition theorems of symmetric products 
and other functors of complexes", Ann. Math. {\bf 68} (1958), 54--80.
\bibitem{}[GKY] M.A. Guest, A. Koslowski, K. Yamaguchi, ``The space of
polynomials with roots of bounded multiplicity'', preprint  math.AT/9807053
(archives).
\bibitem{} [K], S. Kallel, ``Particle spaces and generalized Poincar\'e
Dualities'', preprint 98, math.AT/9810067 (archives).
\bibitem{} [Kr] D. Karagueuzian, Thesis, Stanford University 1994.
\bibitem{} [KP] D. Kahn, S. Priddy, ``On the homology of non-connected
monoids and their associated groups'', Comment. Math. Helv. {\bf 47}
(1972) 1--14.
\bibitem{} [L1] K. Lesh, ``Infinite loop spaces from group theory'',
Math. Z. {\bf 225} (1997), no.3, 467--483.
\bibitem{} [L2] K. Lesh, ``Identification of infinite loop spaces arising
from group theory'', preprint 97.
\bibitem{} [MS] D. McDuff, G. Segal, ``Homology fibrations and the group
completion theorem'', Invent. Math. {\bf 31} (1976), 279--284.
\bibitem{} [Mi] R.J. Milgram, ``Iterated loop spaces'',
Annals of Math.(2) {\bf 84} 1966, 386--403.
\bibitem{}[Ma] J.P. May, ``The geometry of iterated loop spaces'', Springer lecture notes, {\bf 271} (1972).
\bibitem{}[S1] G. Segal, ``The topology of spaces of rational functions",
Acta. Math., {\bf 143}(1979), 39--72.
\bibitem{}[S2] G. Segal, ``Configuration spaces and iterated loop spaces",
Invent. Math., {\bf 21}(1973), 213--221.

\end{thebibliography}
\bibliographystyle{plain}

\flushleft{
Sadok Kallel\\
Dept. of Math., \#121-1884 Mathematics Road\linebreak
U. of British Columbia, Vancouver V6T 1Z2 \linebreak
{\sc Email:} skallel@math.ubc.ca}

\end{document}